\documentclass[11pt]{amsart}

\usepackage[utf8]{inputenc}
\usepackage[english]{babel}
\usepackage[foot]{amsaddr}
\usepackage{hyperref}
\usepackage{longtable}
\usepackage{multirow}
\input xy
\xyoption{all}
\usepackage{mathrsfs}
\usepackage[dvipsnames]{xcolor}
\usepackage[a4paper]{geometry}
\usepackage{indentfirst}

\usepackage[pdftex]{graphicx}

\usepackage{amsmath}
\allowdisplaybreaks  % para que pueda partir fórmulas que ocupan más de una línea, necesita el paquete anterior
\usepackage{amssymb} % para cargar algunos símbolos como \blacksquare y \square
\usepackage{amsfonts} % para cargar algunas fuentes en estilo matemático
\usepackage{enumerate} % Teoremas (se pueden definir todos los que se necesiten):

\numberwithin{equation}{section}

\newtheorem{theorem}{Theorem}[section]
\newtheorem{proposition}[theorem]{Proposition}
\newtheorem{definition}[theorem]{Definition}
\newtheorem{lemma}[theorem]{Lemma}
\newtheorem{corollary}[theorem]{Corollary}

\newtheorem{remark}[theorem]{Remark}

\newcommand{\R}{\mathbb{R}}

\newcommand{\N}{\mathbb{N}}

\newcommand{\Lip}{\operatorname{Lip}}

\newcommand{\LIP}{\operatorname{LIP}}
\newcommand{\md}{\operatorname{md}}
\newcommand{\ud}{\operatorname{d}}
 \newcommand{\Leb}[1]{{\mathscr L}^{#1}}

\providecommand*{\vint}[1]{\mathchoice
	{\mathop{\vrule width 5pt height 3 pt depth -2.5pt
			\kern -9pt \kern 1pt\intop}\nolimits_{\kern -5pt{#1}}}
	{\mathop{\vrule width 5pt height 3 pt depth -2.6pt
			\kern -6pt \intop}\nolimits_{\kern -3pt{#1}}}
	{\mathop{\vrule width 5pt height 3 pt depth -2.6pt
			\kern -6pt \intop}\nolimits_{\kern -3pt{#1}}}
	{\mathop{\vrule width 5pt height 3 pt depth -2.6pt
			\kern -6pt \intop}\nolimits_{\kern -3pt{#1}}}}

\makeatletter
\@namedef{subjclassname@2020}{%
  \textup{2020} Mathematics Subject Classification}
\makeatother

\title{Stepanov Theorem for mappings between metric spaces.}

\author[I. Caama\~{n}o]{Iv\'an Caama\~{n}o}
\address{Institute of Mathematics of the Polish Academy of Sciences, Warsaw, Poland.}
\email{icaamanoaldemunde@impan.pl}

\thanks{The research is partially supported by grant PID2022-138758NB-I00 (Spain) and the Institute of Mathematics of the Polish Academy of Sciences.}


\begin{document}

\maketitle
\begin{abstract}
By combining the ideas of Kirchheim \cite{Kirchheim} and Cheeger \cite{che99} a notion of metric differentiability for mappings between metric spaces was studied in \cite{CDJPS}. A Rademacher-type theorem was introduced there, and this paper will study the validity of a Stepanov-type theorem under the same assumptions.
\end{abstract}

\section{Introduction}

The well-known theorem of Rademacher was generalized to the setting of mappings into a metric space through the notion of \emph{metric differentiability}, introduced by Kirchheim \cite{Kirchheim}. More specifically, Kirchheim proved in \cite{Kirchheim} that every Lipschitz map $f:E\to X$ from a measurable set $E\subset \R^n$ into a metric space is \emph{metrically differentiable}: for  $\Leb{n}$-a.e. $x\in E$ there exists a seminorm $\md_xf$ on $\R^n$ so that
\begin{align}\label{eq:kirch-chart}
	d(f(y),f(z))=\md_xf(y-z)+o(d(y,x)+d(x,z)).
\end{align}
As an illustration of their use, metric differentials give rise to area and co-area formulae, and can be used to show that the Lipschitz maps in the definition of a rectifiable set can be decomposed into bi-Lipschitz maps on a finer atlas, see \cite{Kirchheim,amb-kir00}. 

Parallel developments in analysis on metric spaces lead to a generalized Rademacher theorem for Lipschitz mappings in metric measure spaces $(X,d,\mu )$ with $\mu$ doubling and supporting a Poincaré inequality. This  approach was first studied by Cheeger in \cite{che99}, introducing what have come to be known as Lipschitz differentiability spaces (LDS for short) in recent literature. These are spaces covered by countably many \emph{Cheeger charts}. A Cheeger chart is a pair $(U,\varphi)$ consisting of a Borel set $U\subset X$ with $\mu(U)>0$ and $\varphi\in \LIP(X,\R^n)$ such that every $f\in \LIP(X)$ is differentiable $\mu$-a.e. on $U$ with respect to $(U,\varphi)$: for $\mu$-a.e. $x\in U$ there exists a unique linear map $\ud_xf\in (\R^n)^*$ so that
\begin{align}\label{eq:cheeger-chart}
	f(y)-f(x)=\ud_xf(\varphi(y)-\varphi(x))+o(d(x,y)).
\end{align}
Cheeger showed that spaces endowed with a doubling measure supporting some Poincar\'{e} inequality (in short, PI-spaces) admit a countable covering by such charts, initially referred to as a measurable differentiable structure in early literature (see e.g. \cite{Keith,Keith04}. Recent works have aimed to understand the nature of these spaces, even without assuming $\mu$ doubling and a Poincaré inequality on $X$ (see e.g. \cite{Bate,BaLi,BaSp,IPS}).

Recall that an $n$-dimensional chart $(U,\varphi )$ in $X$ consists on a Borel set $U\subset X$ and a Lipschitz map $f:X\rightarrow \R^n$. In \cite{CDJPS} the notions of Cheeger's differentiability and Kirchheim's metric differentiability were combined to study regularity of Lipschitz mappings $f:X\rightarrow Y$ whenever $(X,d,\mu )$ is a metric measure space and $(Y,d_Y)$ a metric space, introducing the notion of metric differentials with respect to a chart.
\begin{definition}\label{def:metricdifferentiability}
	Given a map $f:X\to Y$ into a metric space $Y$, we say that $f$ admits a (weak) metric differential with respect to the chart $(U,\varphi)$ if there exists a Borel map $\md f: U\to \mathrm{sn}^n$ satisfying, for $\mu$-a.e. $x\in U$,
	\begin{equation}\label{eq:kir-chart-diff}
		\limsup_{U\ni y\to x}\frac{|d_Y(f(y),f(x))-\md_xf(\varphi(y)-\varphi(x))|}{d(x,y)}=0. %,\quad y,z\in X.
	\end{equation}
	Here $\mathrm{sn}^n$ denotes the collection of seminorms in $R^n$.
\end{definition}
See also \cite{CheKleiSc} for an alternative approach to metric differentiation of mappings between metric spaces which relies on metric differentiation along curves, which is compatible with definition \ref{def:metricdifferentiability} whenever both exist, see \cite[Remark 3.3]{CDJPS}. A similar approach to Definition \ref{def:metricdifferentiability} was studied in \cite{GT} for maps defined on strongly rectifiable metric spaces. 
 
 The study of metric differentials has become an indispensable tool in studying rectifiability of metric spaces, and in \cite{CDJPS} such connection was studied for the notion of metric differentiability in Definition \ref{def:metricdifferentiability}. Recall that a metric measure space $X=(X,d,\mu)$ is called $n$-rectifiable if $\mu\ll\mathcal H^n$ and $\mu\left( X\backslash\bigcup_{i=1}^\infty \psi_i(E_i) \right) =0$
 for some countable family of Lipschitz maps $\psi_i\in \LIP(E_i,X)$ defined on $\Leb{n}$-measurable sets $E_i\subset \R^n$.
 
The classical result of Rademacher in the euclidean setting was strengthened by Stepanov in \cite{stepanov}, where differentiability almost everywhere was proved for points in the set
$$S(f):=\{ x\in\R^n:\limsup_{y\to x}\frac{|f(x)-f(y)|}{|x-y|}<\infty \},$$
where $f$ may not be Lipschitz. The results in \cite{MalZaj} show that in more general settings one can obtain an analogue of Stepanov's theorem whenever a Rademacher type theorem holds. Metric differentiability almost everywhere in the sense of Definition \ref{def:metricdifferentiability} of Lipschitz maps into a metric space was proved in \cite{CDJPS} under the (necessary) assumption that $X$ is a metric measure space admitting a rectifiable decomposition, i.e., $X$ can be expressed, up to a null set, as a union of rectifiable charts. Recall that a chart $(U,\varphi )$ is $n$-rectifiable if $\varphi :U\rightarrow \varphi (U)\subset \R^k$ is bi-Lipschitz and $\varphi_\# (\mu_{| U}) \ll \mathcal{L}^{k}_{|\varphi (U)}$. Thus one may ask if the generalization to the set $S(f)$ still holds for non-Lipschitz maps in this setting, and in this work we obtain such result:

\begin{theorem}\label{thm:stepanoff}
	Let $(X,d,\mu )$ be a metric measure space whose porous sets are null and $(U,\varphi )$ a rectifiable chart in $X$. Then any mapping $f:X\rightarrow Y$, where $(Y,d_Y)$ is a metric space, is metrically differentiable with respect to $(U,\varphi )$ almost everywhere on $U\cap S(f)$.
\end{theorem}

The paper is organized as follows: In section 2 we establish basic definitions and notation and give a brief introduction on Lipschitz differentiability spaces, while Section 3 is devoted to the notion of rectifiability. In Section 4 we recall the definition of metric differentiability mentioned above, and study some of its properties in order to prove the main result in Section 5.

\section{Lipschitz differentiability spaces.}
\subsubsection*{Background.} 
Throughout this paper the triplet $(X,d,\mu )$ denotes a complete separable metric space endowed with a measure $\mu$ which is Borel regular and finite on bounded sets. In particular, $\mu$ is Radon (see \cite[Corollary 3.3.47]{HKST}).

Given two measures $\mu$ and $\nu$, we say that $\mu$ is {\em absolutely continuous} with respect to $\nu$ (denoted $\mu\ll\nu$) if whenever $A\subset X$ such that $\nu (A)=0$ then $\mu (A)=0$.

For a mapping $f:(X,d_X)\rightarrow (Y,d_Y)$ between metric spaces we define the pointwise Lipschitz constant as
$$\Lip f(x):=\limsup_{y\to x}\frac{d_Y(f(x),f(y))}{d_X(x,y)}.$$
and $\Lip  f(x)=0$ when $x$ is an isolated point. Denote the set of points $X$ where $f$ is pointwise Lipschitz by
$$S(f):=\{ x\in X : \mathrm{Lip}\, f(x)<\infty \}$$
and let $\mathrm{LIP}(X,Y)$ be the space of Lipschitz mappings. For the particular case of $Y=\R$, we use the notation $\mathrm{LIP}(X)$.

We recall here an embedding result and a Lipschitz extension result that will be useful later on.

\begin{lemma}[Kuratowski embedding Theorem, \cite{Kuratowski}]\label{Kuratowski} Any metric space $(Y,d_Y)$ admits an isometric embedding into the Banach space $\ell^\infty (Y)$.
\end{lemma}

\begin{lemma}\textbf{\em \cite[Corollary 4.1.7]{HKST}} \label{extension}
	Let $(X,d)$ a metric space. Given a Lipschitz mapping $f:U\subset X\rightarrow \ell^\infty (Y)$, where $Y$ is any set, then there exists a Lipschitz map $\hat{f}:X\rightarrow \ell^\infty (Y) $ such that $f(x)=\hat{f}(x)$ for every $x\in U$.
\end{lemma}

Some of the classical results in analysis that do not involve derivatives were extended in the late 1970s to metric spaces endowed with a doubling measure (called in \cite{CW} spaces of homogeneous type). We refer to \cite{Hbook} to an extensive introduction to this topic. In particular, one of the important results developed in this setting is the Lebesgue differentiation theorem.
\begin{theorem}\label{thm:lebesgue}\textbf{\em \cite[Theorem 1.8]{Hbook}}
	Let $(X,d,\mu )$ be a metric measure space with $\mu$ doubling,and $f\in L^1_{loc}(X)$. Then
	$$\lim_{r\to 0}\vint{B(x,r)}\vert f(y)-f(x)\vert \; d\mu (y) =0$$
	for almost every $x\in X$. In particular
	$$\lim_{r\to 0}\vint{B(x,r)}f\; d\mu =f(x)$$
	for almost every $x\in X$. 
\end{theorem}
 We will not make use of the doubling assumption in this work, but some of its consequences also follow from a weaker notion, the infinitessimally doubling condition, which will be present as an underlying consequence of the geometry required throughout the upcoming sections.

\begin{definition}
	The measure $\mu $ is \textit{infinitesimally doubling at $x\in X$} if
	$$\limsup_{r\to 0}\frac{\mu (B(x,2r))}{\mu (B(x,r))}<\infty .$$
	We say that $\mu$ is an \textit{infinitesimally foubling measure} on $X$ if it is infinitesimally doubling for almost every $x\in X$.
\end{definition}
Under the assumption that $\mu$ is infinitesimally doubling at almost every point, Lebesgue differentiation Theorem holds, see for example, in \cite[Theorem 3.4.3 and page 81]{HKST}. In particular, we are interested in the geometric consequence that almost every point in a Borel set is a density point.

\begin{definition}\em \label{def:AppCont}
		Let $(X,d,\mu )$ be a metric measure space and $(Y,d_Y)$ be a metric space. We define the \textit{density} of a set $E\subset X$ at a point $x\in E$ as
	$$\limsup_{r\to 0^+}\frac{\mu (B(x,r)\cap E)}{\mu (B(x,r))},$$
	and we say that $x$ is a \textit{density point} of $E$ if such limit equals $1$.
	 We say that $z_0\in Y$ is the\textit{ approximate limit }of a mapping $f:X\rightarrow Y$ at a point $x\in X$, and denote it by 
	$$z_0=\underset{y\to x}{\mathrm{ap}\lim}f(y),$$
	if for every $\varepsilon >0$ the set $\{ y\in X:d_Y(f(y),z_0)\ge \varepsilon \}$ has density zero at $x$.
\end{definition}
As mentioned above, if $X$ is such that Theorem \ref{thm:lebesgue} holds, then almost every point of a Borel set in $X$ is a density point, and this is then the case whenever $\mu$ is infinitesimally doubling almost everywhere. Although this restriction on $\mu$ will not be a standing assumption, we will notice later on that it will be a consequence of another geometric restriction which is closely related to Cheeger differentiability, see Remark \ref{rem:porousdoub}.

\subsubsection*{Lipschitz charts and weak Lipschitz charts.}
This section will introduce some background on Lipschitz differentiability spaces, and for that we first present the various notions of charts that we need to consider in order to establish a structure on a metric measure space, enabling us to accurately define a differential.

\begin{definition}[Lipschitz charts]
	We say that $(U,\varphi )$ is a \emph{Lipschitz chart} in $X$ if $U\subset X$ is Borel and $\varphi :X\rightarrow \R^n$  is a Lipschitz mapping such that for every function $f\in \mathrm{LIP}(X)$ and almost every $x\in U$ there exists an unique linear map $d_xf:\R^n\rightarrow \R$ such that
	$$\mathrm{Lip} (f-d_xf\circ \varphi )\, (x)=0,$$
	that is,
	\begin{equation}\label{eq:cheegerchart}
		\limsup_{y\to x}\frac{\vert f(x)-f(y)-d_xf(\varphi (x)-\varphi (y))\vert}{d(x,y)}=0.
	\end{equation}
	We call the number $n$ the dimension of the chart.
\end{definition}
The notion of Lipschitz chart is often referred to in the literature as a \textit{Cheeger chart}, since it was first considered in \cite{che99}. 

\begin{definition}[Weak Lipschitz charts]
	We say that $(U,\varphi )$ is a \emph{weak Lipschitz chart} in $X$ if $U\subset X$ is Borel and $\varphi  :X\rightarrow \R^n$  is a Lipschitz mapping such that for every function $f\in \mathrm{LIP}(X)$ and almost every $x\in U$ there exists an unique linear map $d_xf:\R^n\rightarrow \R$ such that
	$$\mathrm{Lip} (f-d_xf\circ \varphi )_{|U }\, (x)=0,$$
	that is,
	\begin{equation}\label{eq:weakcheegerchart}
		\limsup_{U\ni y\to x}\frac{\vert f(x)-f(y)-d_xf(\varphi (x)-\varphi (y))\vert}{d(x,y)}=0.
	\end{equation}
\end{definition}
In the literature, another notion is often considered instead of that of a weak Lipschitz chart, and that is an \textit{approximate Lipschitz chart}, where for every Lipschitz function $f\in \mathrm{LIP}(X)$ and a.e. $x\in U$ there exists an unique linear map $d_xf:\R^n\rightarrow \R$ so that
\begin{equation}\label{eq:apcheegerchart}
	\underset{y\to x}{\mathrm{ap}\lim}\frac{\vert f(x)-f(y)-d_xf(\varphi (x)-\varphi (y))\vert}{d(x,y)}=0
\end{equation}
We will not employ this notion here, but some of the references we provide may utilize it, so we provide an easy argument that will allow us to make use of the results present in the literature that are stated for approximate Lipschitz charts.

Notice that, from Definition \ref{def:AppCont}, one can see that $L$ is the approximate limit of $f:X\rightarrow \R$ at a point $x\in X$ if there exists a Borel set $A\subset X$ such that $x$ is a density point of $A$ and 
$$\lim_{A\ni y\to x}f(y)=L.$$
This implies that, whenever $(X,d,\mu)$ satisfies the condition that almost every point of a Borel subset is a density point, then for a chart $(U,\varphi)$ in $X$, choosing $A=U$ as above, condition \eqref{eq:weakcheegerchart} yields \eqref{eq:apcheegerchart} for almost every $x\in U$. Consequently, the notion of a weak Lipschitz chart is stronger than that of an approximate Lipschitz chart. It is worth noting that the assumption of density points representing almost every point for any Borel subset of $X$ is needed here and will be addressed in the subsequent discussion concerning porous sets. We will observe there that this assumption holds true whenever porous sets are null. Therefore, references utilizing the concept of an approximate Lipschitz chart can be applied whenever porous sets are null.

\begin{definition}\em \label{def:lds}
	We say that a metric measure space $(X,d,\mu )$ is a \emph{Lipschitz differentiability space}, LDS in short, (resp. weak LDS) if there exists a countable family $\{ (U_i ,\varphi_i) \}_{i\in\N}$ of Lipschitz charts (resp. weak Lipschitz charts) of arbitrary dimension and such that $\mu\left( X\backslash \bigcup_i  U_i\right) =0$.
\end{definition}
Given a $n$-dimensional Lipschitz chart $(U_i ,\varphi_i )$ of an LDS, a result from De Philippis, Marchese and Rindler \cite[Theorem 4.1.1]{DMR} concludes that
$$(\varphi_i )_\# (\mu_{| U_i}) \ll \mathcal{L}^n_{|\varphi _i (U_i )}.$$
A related result was provided by Kell and Mondino in \cite[Theorem 1.3]{KeMo}. Namely, if in addition to $X$ being an LDS, it also admits an atlas of \textit{bi-Lipschitz charts}, that is, a chart $(U_i,\varphi_i)$ where $\varphi_i:U_i\rightarrow \R^{n_i}$ is bi-Lipschitz,
then there exists a collection of Borel sets $\{ E_j\}_{j\in\N}$ covering $X$ up to a measure zero set, bi-Lipschitz equivalent to Borel sets in $\R^{N_j}$ and
$$\mu_{| E_j} \ll \mathcal{H}^{n_j}.$$
\subsubsection*{Porous sets.}
It immediately follows that an LDS is, in particular, a weak LDS. However, the converse might fail. An additional assumption was considered in [16] to characterize when the notions of weak LDS and LDS coincide. Here we briefly go through the ideas presented in that paper in order to provide some of the results that will be used later on.
\begin{definition}\em
	Let $(X,d,\mu )$ be a metric measure space. Given $S\subset X$ and  $x\in S$, $S$ is called $\eta$-\textit{porous} at $x$, for $\eta >0$, if there exists a sequence $x_j\to x$ with 
	$$d(x_j,S):=\inf \{ d(x_j,y):y\in S\} \ge \eta d(x_j,x)\quad\mbox{ for all }j\in\N.$$ 
	We say that $S$ is \textit{porous} if every $x\in S$ is $\eta$-porous for some $\eta>0$. 
\end{definition}
In \cite[Section 2]{BaSp} or \cite[Section 4]{Bate} it was shown that, given a porous set $S\subset X$, one can construct a Lipschitz function in $X$ that is almost nowhere differentiable in $S$, and thus any porous set must have measure zero whenever $X$ is a Lipschitz differentiability space, see \cite[Theorem 2.4]{BaSp}.  

\begin{remark}\label{rem:porousdoub}
	The condition of porous sets having measure zero is intermediate between the doubling and the pointwise doubling conditions, that is, if the measure is doubling, then all porous sets are null, as mentioned in \cite[Remark 2.9]{GuyKl}. Moreover, according to \cite[Theorem 3.6 (iv)]{MMPZ}, if all porous sets in $X$ have measure zero then $\mu$ is pointwise doubling almost everywhere.
\end{remark}
It also holds that porous sets being null is a sufficient condition to self-improve a weak-Lipschitz chart into a Lipschitz chart (see \cite[Lemma 2.6]{BaLi18}, or \cite[Proposition 2.8]{BaSp}; in this last reference, the authors prove self-improvement of an approximate Lipschitz chart into a Lipschitz chart). To show this, we first give the following lemma.
\begin{lemma}\label{lem:denseyprime}
	Let $(X,d,\mu )$ be a metric measure space such that all porous sets have zero measure and $U\subset X$ Borel. Then for almost all $x\in U$ and for every $\varepsilon >0$ there exists $r>0$ such that for $y\in B(x,r)$ there exists $z(y)\in U$ with $d(y,z(y))\leq \varepsilon d(x,y)$.
\end{lemma}
A proof of the above claim can be found in \cite[Theorem 3.5]{Keith} or \cite[Proposition 2.9]{BRZ} when the measure is assumed to be doubling. Fortunately, the same ideas can be adapted for non-doubling measures, as long as porous sets have measure zero.

{\em\noindent Proof of Lemma \ref{lem:denseyprime}}\quad
Let $x\in U$ be a density point. Assume by contradiction that there exists $\varepsilon >0$ and a sequence $y_k\in X$ converging to $x$ such that $B(y_k,\varepsilon d(x,y_k))\cap U)=\emptyset$ for all $k\in\N$. We can assume that $\varepsilon \leq 1$ and then by the density condition of $x$ where the radius $r_k\to 0$ is chosen as the sequence $r_k=(1+\varepsilon )d(x,y_k)$
\begin{eqnarray}\label{eq:lemmadenseyprime}
	0&=& \limsup_{r\to 0}\frac{\mu (B(x,r)\setminus U)}{\mu (B(x,r))}=\limsup_{k\to \infty}\frac{\mu (B(x,(1+\varepsilon )d(x,y_k))\setminus U)}{\mu (B(x,(1+\varepsilon )d(x,y_k)))}\nonumber \\
	&\geq & \limsup_{k\to \infty}\frac{\mu (B(y_k,\varepsilon d(x,y_k))\setminus U)}{\mu (B(x,(1+\varepsilon )d(x,y_k)))}
	= \limsup_{k\to \infty}\frac{\mu (B(y_k,\varepsilon d(x,y_k)))}{\mu (B(x,(1+\varepsilon )d(x,y_k)))}\nonumber \\
	&\geq & \limsup_{k\to \infty}\frac{\mu (B(y_k,\varepsilon d(x,y_k)))}{\mu (B(x,2d(x,y_k)))}.
\end{eqnarray}
On the other hand, as explained in the proof of \cite[Proposition 2.8]{BaSp}, using \cite[Theorem 3.6]{MMPZ} we have that, for almost every $x\in U$ and for all $\varepsilon >0$
$$\liminf_{y\to x}\frac{\mu (B(y,\varepsilon d(x,y))}{\mu (B(x,2d(x,y))}>0,$$
which leads to a contradiction with \eqref{eq:lemmadenseyprime} for almost all density points of $U$. Finally, the claim follows since the hypothesis of porous sets being null implies that $\mu$ is pointwise doubling almost everywhere (Remark \ref{rem:porousdoub}), and thus almost all points in $U$ are density points by \cite[Theorem 3.4.3 and p. 81]{HKST}.
\hfill$\blacksquare$

We present now a direct conclusion from Lemma \ref{lem:denseyprime} that adds comfort to the technicalities involving this result. First, as in the hypotheses of Lemma \ref{lem:denseyprime} assume that porous sets in $X$ are null and let $U\subset X$ Borel.  Given $x\in U$ that satisfies the conclusion of Lemma \ref{lem:denseyprime} and a sequence $\{ y_k\}_{k\in\N}\subset X$ converging to $x$, there exists a sequence $\varepsilon_k>0$ converging to $0$ such that $d(x,y_k)<\varepsilon_k$ for all $k\in\N$. Now applying Lemma \ref{lem:denseyprime} to $x$ and $\varepsilon_k$ for each, $k$ there exist radii $\{ r_k\}_{k\in\N}$ such that for each $y\in B(x,r_k)$ there exists $z(y)\in U$ with $d(y,z(y))<\varepsilon d(x,y)$. Going back to the sequence $\{ y_k\}_{k\in\N}$, we can choose a subsequence $\{ y_{k_j}\}_{j\in\N}$ such that $y_{k_j}\in B(x,r_j)$ for all $j\in\N$. Then, denoting $z_j:=z(y_{k_j})$ we have 
\begin{equation}
	\label{eq:porouslimit}
	\limsup_{j\to \infty }\frac{d(y_{k_j},z_j )}{d(x,y_{k_j})}=0.
\end{equation}
This yields that $\mathrm{Lip}\, f_{|U}=\mathrm{Lip}\, f$ almost everywhere in $U$ for any Lipschitz mapping $f:X\rightarrow \R$. Indeed, choosing $z_j\in U$ as in \eqref{eq:porouslimit}
\begin{eqnarray*}
	\mathrm{Lip}\, f (x)&=&\limsup_{j\to \infty }\frac{|f(x)-f(y_{k_j})|}{d(x,y_{k_j})} 
	\\
	&\leq &\limsup_{j\to \infty}\left(\frac{|f(x)-f(z_j )|}{d(x,z_j )}\frac{d(x,z_j )}{d(x,y_{k_j})}
	+ \frac{|f(y_{k_j})-f(z_j )|}{d(x,y_{k_j})} \right)\\
	&\leq &
	\limsup_{j\to \infty}\left(\frac{|f(x)-f(z_j )|}{d(x,z_j )}
	+ \mathrm{LIP}(f)\frac{d(y_{k_j},z_j )}{d(x,y_{k_j})} \right)
	\\
	&=&\mathrm{Lip}\, f_{|U}(x).
\end{eqnarray*}
Here we used \eqref{eq:porouslimit} together with the triangle inequality to obtain $\limsup_{j\to \infty}\frac{d(x,z_j )}{d(x,y_{k_j})}\leq 1$. Since the inequality $\mathrm{Lip}\, f_{|U}\leq \mathrm{Lip}\, f$ holds in general this proves that $\mathrm{Lip}\, f(x)=\mathrm{Lip}\, f_{|U}(x)$ at almost every point $x\in U$. This fact yields the self-improvement of a weak Lipschitz chart into a Lipschitz chart, giving the following characterization.

\begin{corollary}\textbf{\em \cite[Theorem 2.4 and Proposition 2.8]{BaSp}}\label{cor:wlds-lds}
	Let $(X,d,\mu )$ be a metric measure space. Then $X$ is a Lipschitz differentiability space if, and only if, it is a weak Lipschitz differentiability space and porous sets are null.
\end{corollary}

We point out that a union of Lipschitz differentiability spaces may not be a Lipschitz differentiability space, as shown by D. Bate and S. Li in \cite[Page 5]{BaLi} where they consider the set of $\R^2$ given by
$$(\{ 0\}\times [0,1])\cup \{ (x,p/2^n): n\in\N, 1\leq p <2^n\mbox{ odd}, \pm x\in [2^{-n}-4^{-n},2^{-1}]\},$$
equipped with the $\mathcal{H}^1$-measure, which despite being a countable union of intervals, the map $|x|$ is nowhere differentiable at $\{ 0\}\times [0,1]$. In fact, this vertical segment is a porous set of positive $\mathcal H^1$-measure. 

In contrast to the above example, assume that $(X,d,\mu )$ is a metric measure space such that porous sets have measure zero and $X=\bigcup_{i\in\N} X_i$, where $(X_i,d|_{X_i},\mu|_{X_i})$ is an LDS for each $i\in\N$. Then, since unions of weak Lipschitz differentiability spaces are also weak Lipschitz differentiability spaces, Corollary \ref{cor:wlds-lds} implies that $X$ is also an LDS.

\section{Rectifiability.}\label{sec:rectifiability}
In order to give a formal definition of rectifiability, we briefly recall the notion of {\em Hausdorff measure}. For $s> 0$, first fix $\delta >0$ and consider, for any set $E\subset X$,
$$\mathcal{H}_\delta^s(E):=\inf\left\{ \sum_{i=1}^\infty \mathrm{diam}(E_i)^s \right\},$$
where the infimum is taken over all countable covers of $E$ by sets $E_i\subset X$ with $\mathrm{diam}(E_i)<\delta$ for all $i\in\N$. Then the $s$-\textit{dimensional Hausdorff measure} is defined as
$$\mathcal{H}^s(E):=\lim_{\delta\to 0}\mathcal{H}_\delta^s(E).$$

We present two different ways, but very related, of decomposing a metric measure space into a family of sets Lipschitz equivalent to Euclidean sets. This notion is known as rectifiability and we first present the classical definition due to Federer.

\begin{definition}[Rectifiable space] \label{rectispace}
	A metric measure space $(X,d,\mu)$ is {\em $k$-rectifiable} if there exists a countable family of Lipschitz mappings $\psi_i :E_i\subset \R^k\rightarrow X$ defined on measurable sets $E_i$ such that
	$$\mu\left( X\backslash\bigcup_{i=1}^\infty \psi_i(E_i) \right) =0\quad\text{ and }\quad \mu\ll\mathcal{H}^k.$$
\end{definition}

\begin{remark} \label{rem:bilippieces}
	According to a result of Kirchheim \cite[Lemma 4]{Kirchheim}, Lipschitz functions that parametrize a $k$-rectifiable set can be chosen to be bi-Lipschitz, that is, if $X$ is  $k$-rectifiable with decomposition  $\psi_i :E_i\subset \R^k\rightarrow X$, then there exist $\{ E_{i,j}\}_{j\in\N}\subset E_i$ for all $i$ so that
	\begin{enumerate}
		\item[$(i)$] $E_{i,j}\subset E_i$ is Borel for each $j\in\N$ and $\mu\left( \psi_i (E_i)\backslash \underset{j\in\N}{\bigcup}\psi_i (E_{i,j})\right) =0$ for each $i\in\N$ and
		\item[$(ii)$] $\psi_{i,j}:=\psi_{i|_{E_{i,j}}}$ is bi-Lipschitz.
	\end{enumerate}
	In the Euclidean setting, this follows from \cite[Lemma 3.2.2]{Fe}. Observe that a $k-$rectifiable metric measure space has a natural decomposition into charts. Namely, using the notation in Definition \ref{rectispace}, if $(X,d,\mu)$ is $k-$rectifiable, it can be written as a countable collection of $k-$dimensional charts $\{(U_i,\varphi_i)\}_{i\in\N}$ where
	\[
	U_i=\psi_i(E_i)\text{ and } \varphi_i=\psi_i^{-1} \text{ are bi-Lipschitz} .
	\]
	This fact inspires the following definition, see for example \cite[Section 5]{IPS}.
\end{remark}

\begin{definition}[Rectifiable chart] \label{RectiChart}
	We say that a  $k-$dimensional chart $(U,\varphi)$ is {\em $k-$rectifiable} if $\varphi$ is bi-Lipschitz on $U$ and
	\begin{equation}\label{eq:acmeasures}
		\varphi_\# (\mu_{| U}) \ll \mathcal{L}^{k}_{|\varphi (U)}.
	\end{equation}
	We say that a metric measure space $(X,d,\mu )$ \textit{admits a rectifiable decomposition} if there exists a countable family $\{ (U_i ,\varphi_i )\}_{i\in\N}$ of rectifiable charts such that $\mu\left( X\backslash \bigcup_i  U_i\right) =0$.
\end{definition}

\begin{lemma}\label{lem:rectweaklip}
	Let $(X,d,\mu )$ be a metric measure space and let $(U,\varphi )$ be a rectifiable chart of dimension $n$ in $X$. Then there exist a countable collection of Borel sets $V_i\subset U$ with $\mu\left( U\backslash \bigcup_iV_i\right) =0$ and constants $C_i>0$ so that
	\begin{equation}\label{abscontequiv}
		\frac{1}{C_i }\mathcal{L}^{k}_{|\varphi_i  (V_i)}\leq \varphi_i \, _\# (\mu_{| V_i}) \leq C_i \mathcal{L}^{k}_{|\varphi_i  (V_i)}
	\end{equation}
	where $\varphi_i =\varphi_{|_{V_i}}:V_i\rightarrow \R^{k}$.
\end{lemma}
\begin{proof}
	The refinement of the chart is given by
	$$
	V_{i}:=\varphi^{-1}\left(\left\{\frac{1}{i}<\frac{d[\varphi_\# (\mu_{|U})]}{d\mathcal{L}^{k}_{|_{\varphi (U)}}}\leq i \right\} \right).
	$$
	By the Radon-Nykodym Theorem we have
	\[
	\varphi_\# (\mu_{|U})(E)=\int_E \frac{d[\varphi_\# (\mu_{|U})]}{d\mathcal{L}^{k}_{|_{\varphi (U)}}} d\mathcal{L}^{k}_{|_{\varphi (U)}}.
	\]
	Therefore, \eqref{abscontequiv} holds and the set 
	$E=\left\{\frac{d[\varphi_\# (\mu_{|U})]}{d\mathcal{L}^{k}_{|_{\varphi (U)}}}=0\right\}$ satisfies
	\[
	\mu_{|U}\left(\varphi^{-1}(E)\right)=\varphi_\# (\mu_{|U})(E)=0,
	\]
	proving that $\{ V_i\}_{i\in\N}$ covers $U$ up to a set of measure zero.
\end{proof}

Due to Remark \ref{rem:bilippieces} it follows that if $X$ is $k$-rectifiable with maps $\psi_i:E_i\rightarrow X$ and \eqref{eq:acmeasures} holds for $\varphi_i =\psi_i^{-1}$ then $X$ admits a rectifiable decomposition. Conversely, since $\mathcal H^n\sim \mathcal L^n$ and $\mathcal \varphi_\# \mathcal H^n\ll \mathcal H^n$ for Lipschitz $\varphi$, then \eqref{eq:acmeasures} is stronger than the absolute continuity condition  in Definition \ref{rectispace}, so a rectifiable chart is, in particular, a rectifiable space. However, $X$ admitting a rectifiable decomposition does not necessarily imply that it is a rectifiable space since the charts of the decomposition are of arbitrary dimension. Now the next result follows naturally.

\begin{lemma}\label{unionrectifiable}
	Let $(X,d,\mu )$ be a metric measure space. Then $X$ admits a rectifiable decomposition if, and only if, there exists a disjoint family $\{ (U_i,d|_{U_i},\mu|_{U_i})\}_{i\in\N}$ of $k_i$-rectifiable spaces  such that $X=Z\cup\bigcup_i U_i$, with $\mu (Z)=0$, and 
	$$\frac{1}{C_{i,j}}\mathcal{L}^{k_i}_{|_{E_{i,j}}}\leq (\psi^{-1}_{i,j} )_\# (\mu_{|_{\psi_{i,j} (E_{i,j})}}) \leq C_{i,j} \mathcal{L}^{k_i}_{|_{E_{i,j}}}$$
	for some constants $C_{i,j}>0$, where $\{ (E_{i,j},\psi_{i,j})\}_{j=1}^\infty $ are given by the $k_i$-rectifiability of $U_i$.
\end{lemma}

By \cite[Theorem 4.1.1]{DMR} it follows that an LDS such that its chart maps are bi-Lipschitz also admits a rectifiable decomposition, just by considering as the rectifiable charts the ones given by the LDS condition. However, the converse is not true \cite[p. 5]{BaLi}. Instead, what one has due to Lemma \ref{lem:rectweaklip} and \cite[Lemma 4.1]{IPS} is a decomposition into weak Lipschitz charts with bi-Lipschitz maps.
Given the hypothesis that porous sets in $X$ are null, we know that weak Lipschitz charts self-improve into Lipschitz charts (see Corollary \ref{cor:wlds-lds}). Therefore, we can deduce the following direct consequence for Lemma \ref{lem:rectweaklip}.

\begin{corollary}\label{cor:finiteLDSequivalences}
	Let $(X,d,\mu )$ be a metric measure space such that all porous sets are null and consider a chart $(U,\varphi )$ of dimension $n$ in $X$. Then $(U,\varphi )$ is a rectifiable chart if and only if it admits a decomposition into Lipschitz charts with bi-Lipschitz chart maps.
\end{corollary}
\begin{corollary}\label{cor:LDSequivalences}
	Let $(X,d,\mu )$ be a metric measure space. Then the following are equivalent:
	\begin{enumerate}
		\item[$(i)$] $X$ admits a rectifiable decomposition and all its porous sets are of measure zero.
		\item[$(ii)$] $X$ is an LDS and admits a rectifiable decomposition.
		\item[$(iii)$] $X$ is an LDS with bi-Lipschitz chart maps.
	\end{enumerate}
\end{corollary}

\section{The metric differential.}\label{sec:metricdiff}
In this section we will follow the ideas of \cite{GT} to define a metric differential for mappings in metric measure space with an atlas, and study the relationship of the generalization of Kirchheim's result \cite[Theorem 2]{Kirchheim} and the rectifiability of a chart.

We recall that a seminorm is a mapping $\mathrm{n}:\R^n\rightarrow \R^+$ which is non-negative, subadditive and absolutely homogeneous. We will denote by $\mathrm{sn}^n$ to the space of all seminorms in $\R^n$ endowed with the following metric
$$D(\mathrm{n}_1,\mathrm{n}_2):=\sup_{\vert v\vert\leq 1}\vert \mathrm{n}_1(v)-\mathrm{n}_2(v)\vert =\mathrm{Lip}(\mathrm{n}_1-\mathrm{n}_2)(0).$$
We also denote
$$||| \mathrm{n} |||:=D(\mathrm n,0)=\mathrm{Lip}(\mathrm{n})(0).$$

\begin{definition}\label{def:md}\em 
	Given a metric measure space $(X,d,\mu )$ and a weak Lipschitz chart $(U,\varphi )$ of dimension $n$ in $X$, we say that a mapping $f:X\rightarrow Y$ with values in a metric space $(Y,d_Y)$ is \textit{weakly metrically differentiable }at a given point $x\in U$ with respect to the chart $(U,\varphi )$ if there exists $\mathrm{md}_x f\in \mathrm{sn}^n$ such that
	\begin{equation}\label{weakmd}
		\limsup_{\underset{y\in U}{y\to x}}\frac{\vert d_Y( f(x),f(y)) - \mathrm{md}_x f (\varphi (x)-\varphi (y))\vert }{d(x,y)}=0.
	\end{equation}
	If $(U,\varphi )$ is a Lipschitz chart then we say that $f$ is \textit{metrically differentiable }at $x\in U$ with respect to the chart $(U,\varphi )$ provided \eqref{weakmd} holds for $y\to x$ with $y\in X$ arbitrary.
\end{definition}
If $X=\R^n$ then $(\R^n,id_{\R^n})$ is a Lipschitz chart in $X$ due to Rademacher's Theorem. Hence, for $f:\R^n\rightarrow Y$, our definition yields the classical definition of metric differentiability.
%\begin{example}

\begin{lemma}\label{lem:mdbounded}
	Let $(X,d,\mu )$ be a metric measure space with $(U ,\varphi )$ a rectifiable chart of dimension $n$ in $X$, and $(Y,d_Y)$ a metric space. Then there exists a constant $C>0$ depending only on $\varphi$ such that, for any mapping $f:X\rightarrow Y$ metrically differentiable at $x\in U$ with respect to $(U ,\varphi )$, one has
	\begin{equation}\label{eq:md=Lip}
		\frac{1}{C}\mathrm{Lip}f(x)\leq ||| \mathrm{md}_xf|||\leq C\mathrm{Lip}f(x).
	\end{equation}
\end{lemma}
\begin{proof}
	Let $x\in U$. As $\varphi$ is Lipschitz there exists $C>0$ such that
	\begin{equation}\label{eqcota1}
		\frac{\mathrm{md}_xf(\varphi (y)-\varphi (x))}{d(x,y)}\leq C\frac{\mathrm{md}_xf(\varphi (y)-\varphi (x))}{\vert \varphi (y)-\varphi (x)\vert }.
	\end{equation}
	On the other hand
	$$\mathrm{Lip}(\mathrm{md}_xf)(0)=\limsup_{z\to 0}\frac{\mathrm{md}_xf(z)}{\vert z\vert }.$$
	Considering $z=\varphi (y)-\varphi (x)$, $y\to x$ implies $z\to 0$ as $\varphi $ is continuous, and then
	$$\limsup_{y\to x}\frac{\mathrm{md}_xf(\varphi (y)-\varphi (x))}{\vert \varphi (y)-\varphi (x)\vert }=\limsup_{z\to 0}\frac{\mathrm{md}_xf(z)}{|z|}=\mathrm{Lip}(\mathrm{md}_xf)(0)=|||\mathrm{md}_xf|||.$$
	This, together with \eqref{eqcota1} and the definition of metric differential at $x$, concludes that
	$$\mathrm{Lip}f(x)=\limsup_{y\to x}\frac{d_Y(f(x),f(y))}{d(x,y)}=\limsup_{y\to x}\frac{\mathrm{md}_xf(\varphi (y)-\varphi (x))}{d(x,y)} \leq C\cdot |||\mathrm{md}_xf|||.$$
	For the first inequality in \eqref{eq:md=Lip} we use the Lipschitz condition of $\varphi^{-1}$ to obtain
	\begin{equation}\label{eqcota2}
		\frac{\mathrm{md}_xf(\varphi (y)-\varphi (x))}{d(x,y)}\geq \frac{1}{C}\frac{\mathrm{md}_xf(\varphi (y)-\varphi (x))}{\vert \varphi (y)-\varphi (x)\vert }.
	\end{equation}
	Then we proceed analogously as we did in the second inequality, but using \eqref{eqcota2} instead of \eqref{eqcota1}.
\end{proof}

Lemma \ref{lem:mdbounded} provides a nice estimate of the metric differential for points in the set $S(f):=\{ x\in X:\mathrm{Lip}\, f(x)<\infty\}$ of a measurable mapping $f:X\rightarrow Y$. Moreover, we can decompose $S(f)$ as the union of the sets
\begin{equation}\label{eq:S(f)decomp}
	E_{k}:=\left\{ x\in S(f): \frac{d_Y(f(x),f(y)) }{d(x,y)}\leq k\text{ if }d(x,y) < \frac{1}{k }\right\} ,
\end{equation}
allowing to locally study $f$ as a Lipschitz map in order to apply a Rademacher-type result to $f_{|E_k}$ and obtain metric differentiability almost everywhere in $E_k$. However, to recover the notion of metric differentiability for $f$ and not just for $f_{|E_k}$ one must go through similar arguments as we did for Corollary \ref{cor:wlds-lds}. We will do so in Theorem \ref{thm:wmd-md} below, but before that we first prove that $E_k$ is measurable for each $k$, in order to have that a.e. point is a density point.

\begin{lemma}\label{lem:stepanoffmeasurable}
	Let $(X,d,\mu )$ be a metric measure space and $(Y, d_Y)$ a metric space. Let $f:X\rightarrow Y$ and $L,\delta >0$. Then the following set is closed:
	$$E:=\left\{ x\in X: \frac{d_Y (f(x),f(y)) }{d(x,y)}\leq L\text{ if }d(x,y)< \delta \right\}$$
\end{lemma}
\begin{proof}
	Consider a sequence $\{ x_n\}_{n=1}^\infty\subset E$ such that $d(x_n,x)\to 0$. Let $y\in X$ be such that $d(x,y)<\delta$ and $0<\varepsilon <\delta -d(x,y)$, then there exists $n_0\in\N$ such that $d(x,x_{n_0})<\varepsilon$. Hence
	$$d(y,x_{n_0})\leq d(x,x_{n_0})+d(x,y)<\varepsilon + d(x,y)<\delta .$$
	Then we have
	\begin{eqnarray*}
		d_Y(f(x),f(y)) &\leq & d_Y(f(x),f(x_{n_0}))+d_Y(f(x_{n_0}),f(y))\leq L(d(x,x_{n_0})+d(x_{n_0},y)) \\
		&\leq & L(2d(x,x_{n_0})+d(x,y))< 2L\varepsilon +Ld(x,y)
	\end{eqnarray*}
	and, as $\varepsilon$ is arbitrarily small, we have $d_Y(f(x),f(y))\leq Ld(x,y)$, proving that $x\in E$, and thus, $E$ is closed.
\end{proof}

\begin{theorem}\label{thm:wmd-md}
	Suppose $(X,d,\mu )$ is a metric measure space such that porous sets are null and let $(U ,\varphi )$ be a rectifiable chart of dimension $n$ in $X$. Let $(Y,d_Y)$ be a metric space and $f:X\rightarrow Y$, then $f$ is metrically differentiable at almost every  point $x\in U\cap S(f)$ such that there exists $\mathrm{md}_x f\in\mathrm{sn}^n$ satisfying
	\begin{equation}\label{eqderivadaU}
		\limsup_{\underset{y\in U}{y\to x}}\frac{\vert d_Y(f(x),f(y))- \mathrm{md}_xf(\varphi (x)-\varphi (y))\vert }{d(x,y)}=0.
	\end{equation}
	In particular, weak metric differentiability a.e. yields metric differentiability a.e. with respect to $(U,\varphi )$.
\end{theorem}
\begin{proof}
	First recall that $S(f)=\bigcup_{k\in\N}E_k$ where $E_k$ are defined as in \eqref{eq:S(f)decomp}, and they are Borel for every $k\in\N$ due to Lemma \ref{lem:stepanoffmeasurable}. 
	
	Fix $k_0\in\N$ and $x\in U\cap E_{k_0}$ satisfying the conclusion of Lemma \ref{lem:denseyprime} for the Borel set $U\cap E_{k_0}$ and such that \eqref{eqderivadaU} holds at $x$.
	Let $\{ y_j\}_{j\in\N}$ be a sequence in $X$ converging to $x$, then by Lemma \ref{lem:denseyprime}, and particularly by the discusion leading to \eqref{eq:porouslimit}, passing to a subsequence if necessary (that we still denote by $y_j$), there exists $\{ z_j\}_{j\in\N}$ in $U\cap E_{k_0}$ such that
	\begin{equation}\label{eq:densepoint}
		\limsup_{j\to \infty }\frac{d(y_j,z_j)}{d(x,y_j)}=0.
	\end{equation}
	In particular, $z_j\to x$. Now for each $j\in\N$ we have
	
	\begin{eqnarray*}
		\frac{\vert d_Y( f(x),f(y_j)) - \mathrm{md}_x f (\varphi (x)-\varphi (y_j))\vert }{d(x,y_j)}
		&\leq & 
		\frac{|d_Y(f(x),f(y_j))-d_Y(f(x),f(z_j))|}{d(x,y_j)}\\
		&+& \frac{\vert d_Y(f(x),f(z_j )) - \mathrm{md}_xf (\varphi (x)-\varphi (z_j))\vert }{d(x,y_j)} \\
		&+& \frac{\vert \mathrm{md}_xf (\varphi (y_j)-\varphi (z_j))\vert }{d(x,y_j)}\\
		&=:& \mathrm{I}_j+\mathrm{II}_j+\mathrm{III}_j
	\end{eqnarray*}
	The proof follows by checking
	$$\limsup_{j\to\infty} \mathrm{I}_j=\limsup_{j\to\infty} \mathrm{II}_j=\limsup_{j\to\infty} \mathrm{III}_j=0.$$
	We now estimate each term separately. \\
	First, recall that $z_j\subset E_{k_0}$, and since $d(y_j,z_j)\to 0$ as $j\to\infty$, for sufficiently large $j$ we can assume $d(y_j,z_j)<1/k_0$. Then, by the definition of $E_{k_0}$ we have
	\begin{eqnarray*}
		\limsup_{j\to \infty }\mathrm{I}_j&:=&\limsup_{j\to\infty } \frac{|d_Y(f(x),f(y_j))-d_Y(f(x),f(z_j))|}{d(x,y_j)}\\
		&\leq& \limsup_{j\to\infty} \frac{d_Y(f(y_j),f(z_j)}{d(x,y_j)}\leq k_0\limsup_{j\to\infty}\frac{d(y_j,z_j)}{d(x,y_j)}=0.
	\end{eqnarray*}
	To estimate $\mathrm{II}_j$ notice that using triangle inequality in \eqref{eq:densepoint} yields
	$$\limsup_{j\to \infty }\frac{d(x,z_j )}{d(x,y_j)}\leq 1,$$
	and then by the hypothesis \eqref{eqderivadaU}, since $z_j\in U$ for every $k\in\N$
	\begin{eqnarray*}
		\limsup_{j\to \infty} \mathrm{II}_j&:=&
		\limsup_{j\to\infty}  \frac{\vert d_Y(f(x),f(z_j ))- \mathrm{md}_xf (\varphi (x)-\varphi (z_j))\vert }{d(x,y_j )}\\
		&=&  \limsup_{j\to\infty}\frac{\vert d_Y(f(x),f(z_j ))- \mathrm{md}_xf (\varphi (x)-\varphi (z_j))\vert }{d(x,z_j )}\frac{d(x,z_j )}{d(x,y_j)} \leq 0.
	\end{eqnarray*} 
	Last, by the Lipschitz condition of $\varphi$, Lemma \ref{lem:mdbounded} applied to $f_{|U\cap S(f)}$ and \eqref{eq:densepoint} we have
	\begin{eqnarray*}
		\limsup_{j\to\infty} \mathrm{III}_j&:=&
		\limsup_{j\to \infty} \frac{\vert \mathrm{md}_xf (\varphi (y_j)-\varphi (z_j))\vert }{d(x,_k)}\\
		&\leq &\limsup_{j\to \infty}\frac{\mathrm{LIP}(\varphi ) ||| \mathrm{md}_xf|||  d(y_j,z_j )}{d(x,y_j)}\\
		&\leq &\limsup_{j\to \infty}\frac{C\mathrm{LIP}(\varphi ) \mathrm{Lip}f(x)d(y_j,z_j)}{d(x,y_j)}=0.
	\end{eqnarray*}
	We conclude  that $f$ is metrically differentiable at $x$. Now for each $k\in\N$ we denote by $N_k$ the set of all points in $U\cap E_k$ such that the conclusion of Lemma \ref{lem:denseyprime} does not hold. Then for each $k\in\N$, since porous sets in $X$ are null and $U\cap E_k$ is measurable by Lemma \ref{lem:stepanoffmeasurable}, $\mu (N_k)=0$. Thus $\mu (S(f)\setminus N)=0$, where $N=\bigcup_k N_k$. We have then proved that, for almost every point $x\in S(f)\cap U$ such that \eqref{eqderivadaU} holds, $f$ is metrically differentiable at $x$.
\end{proof}

\section{Stepanov Theorem.}\label{sec:stepanoff}
We now have all the necessary tools to apply the standard arguments that improve a Rademacher Theorem into a Stepanov Theorem. It was shown in \cite[Section 4]{CDJPS} that metric differentiability almost everywhere in a cheeger chart for Lipschitz maps into any metric space is equivalent to the rectifiability of the chart, making the later a natural assumption for this theory. For the sake of completion, we prove here the easier implication, that is, that a Rademacher theorem holds under the assumption of rectifiability, that will be needed for the proof of \ref{thm:stepanoff}.

\begin{proposition}[\cite{CDJPS}, Proposition 4.2]\label{weakKirchheim}
	Let $(X,d,\mu)$ be a metric measure space and let $(U,\varphi)$ be a $n-$rectifiable chart. Then every Lipschitz mapping $f:X\rightarrow Y$ into a metric space $(Y,d_Y)$ is weakly metrically differentiable almost everywhere with respect to $(U,\varphi )$.
\end{proposition}
\begin{proof}
	Let $g=f\circ \varphi^{-1} :\varphi(U)\rightarrow Y$. Notice that $g$ is a composition of Lipschitz mappings, so it is also Lipschitz. By Kirchheim's Theorem \cite[Theorem $2$]{Kirchheim}, for $\mathcal{H}^n$-almost every $z\in \varphi(U)$ there exists a unique seminorm $ \mathrm{md}_zg$ on $\R^{n}$ such that
	\begin{equation}\label{medieq}
		\lim_{\stackrel{y\to z}{y\in \varphi(U) }}\frac{|d_Y(g(z),g(y))-\mathrm{md}_zg(y-z)|}{|y-z|}=0.
	\end{equation}
	On the other hand, $g(\varphi (x))=f(x)$ for each $x\in
	U$. Fix $x_0\in U$ such that for $z_0:= \varphi(x_0)$ there exists a unique seminorm $ \mathrm{md}_{z_0}g$ on $\R^{n}$ such that \eqref{medieq} holds. As $\varphi$ is continuous, if $x\in U$ and $x\to x_0$, then $\varphi(x)\to z_0$. Therefore
	\begin{eqnarray*}
		&&\lim_{\underset{ x\in U }{x\to x_0 }}\frac{|d_Y(f(x_0),f(x))-\mathrm{md}_zg(\varphi(x)-\varphi(x_0))|}{d(x,x_0)} \\
		&\leq &\lim_{\underset{z=\varphi(x)}{z\to z_0}}\frac{\vert d_Y(g(z_0),g(z))-\mathrm{md}_{z_0}g(z-z_0)\vert}{\frac{1}{C}|z-z_0|}=0.
	\end{eqnarray*}
	where  $C $ is the Lipschitz constant of $\varphi$. Then $\mathrm{md}_{x_0}f:=\mathrm{md}_{\varphi(x_0)}g$ is the metric differential of $f$ at  $x_0\in U$. We finish the proof by noticing that, because $\varphi_\# (\mu_{| U}) \ll \mathcal{L}^{n}_{|\varphi (U)}
	$ and
	\[
	\mathcal{L}^{n}(\{z_0\in \varphi(U): \mathrm{md}_{z_0}g \text{ does not exist}\})=0,
	\]
	we conclude that
	\[
	\mu(\{x_0\in U: \mathrm{md}_{x_0}f \text{ does not exist}\})=0.
	\]
\end{proof}
\begin{theorem}
	Let $(X,d,\mu )$ be a metric measure space whose porous sets are null and $(U,\varphi )$ a rectifiable chart in $X$. Then any mapping $f:X\rightarrow Y$, where $(Y,d_Y)$ is a metric space, is metrically differentiable with respect to $(U,\varphi )$ almost everywhere on $U\cap S(f)$.
\end{theorem}
\begin{proof}
	By Kuratowski's embedding we can suppose without loss of generality that $f:X\rightarrow \ell^\infty (Y)$, as differentiability is invariant under isometries. 
	The goal is to prove that $f$ is metrically differentiable at almost every point of $M :=U\cap S(f)$.
	
	For each $k\in\N$ consider the set
	$$E_{k}:=\left\{ x\in S(f): \frac{\Vert f(x)-f(y)\Vert_{\ell^\infty} }{d(x,y)}\leq k\text{ if }d(x,y) < \frac{1}{k }\right\} ,$$
	and let $\{ B_{k,j}\}_{j=1}^\infty$ be a covering by open sets of $X$ such that $\textrm{diam} (B_{k,j})\leq\frac{1}{k}$. For each $k,j\in\N$ denote $M_{k,j}=M\cap E_k\cap B_{k,j}$. 
	
	Clearly $S(f)=\bigcup_{k} E_k$, and then it is measurable by Lemma \ref{lem:stepanoffmeasurable}. Hence, $M$ is measurable as so it is $U$ by definition, concluding that $M_{k,j}$ is a measurable set for each $k,j\in\N$.
	
	On the other hand $S(f)=\bigcup_{k,j} E_k\cap B_{k,j}$ as $\{ B_{k,j}\}$ is a covering of $X$, and then we have
	$$M =  \bigcup_{k,j=1}^\infty M_{k,j},$$
	so it suffices to prove that, for each $k,j\in\N$, $f$ is metrically differentiable almost everywhere in $M_{k,j}$.
	
	Fix $j,k\in\N$. First, notice that $f|_{M_{k,j}}$ is Lipschitz. Indeed, if $x,y\in M_{k,j}$ then $d(x,y)\leq \frac{1}{k}$ since $x,y\in B_{k,j}$, and then by the definition of $E_k$ one has
	$$\Vert f(x)-f(y)\Vert_{\ell^\infty} \leq kd(x,y).$$
	Then there exists a Lipschitz mapping $\hat{f}:X\rightarrow \ell^\infty (Y)$ such that $\hat{f}|_{M_{k,j}}=f|_{M_{k,j}}$, and by Proposition \ref{weakKirchheim} and Theorem \ref{thm:wmd-md} $\hat{f}$ is metrically differentiable with respect to $(U,\varphi)$ at almost every point $x\in U$.
	
	Let $x\in M_{k,j}$ such that $\hat{f}$ is metrically differentiable at $x$. Then there exists a seminorm $\mathrm{md}_x \hat{f}\in\mathrm{sn}^{\mathrm{dim}(\varphi (U))}$ such that
	$$\lim_{y\to x}\frac{\vert \Vert \hat{f}(x)-\hat{f}(y)\Vert_{\ell^\infty} - \mathrm{md}_x \hat{f} (\varphi (x)-\varphi (y))\vert }{d(x,y)}=0.$$
	If $y\in M_{k,j}$ then $\hat{f}(y)=f(y)$, so we have
	\begin{equation}
		\lim_{\underset{y\in M_{k,j}}{y\to x}}\frac{\vert \Vert f(x)-f(y)\Vert_{\ell^\infty} - \mathrm{md}_x \hat{f} (\varphi (x)-\varphi (y))\vert }{d(x,y)}=0.
	\end{equation}
	We conclude that $f$ is metrically differentiable at $x$ by Theorem \ref{thm:wmd-md}. 
\end{proof}

\end{document}